\documentclass{elsart3-1}

\usepackage{amssymb}
\usepackage{amsmath,color}
\usepackage[english,francais]{babel}

\newtheorem{theorem}{Theorem}[section]
\newtheorem{lemma}[theorem]{Lemma}
\newtheorem{e-proposition}[theorem]{Proposition}

\newtheorem{e-definition}[theorem]{Definition\rm}

\newtheorem{theoreme}{Th\'eor\`eme}[section]

\newtheorem{proposition}[theoreme]{Proposition}

\renewcommand{\theequation}{\arabic{equation}}
\newcommand{\R}{\mathbb{R}}
\newcommand{\E}{\mathrm{E}}
\renewcommand{\P}{\mathrm{P}}

\def\og{\leavevmode\raise.3ex\hbox{$\scriptscriptstyle\langle\!\langle$~}}
\def\fg{\leavevmode\raise.3ex\hbox{~$\!\scriptscriptstyle\,\rangle\!\rangle$}}

\journal{the Acad\'emie des sciences}
\begin{document}
\centerline{}
\begin{frontmatter}


\selectlanguage{english}
\title{LAN property for a linear model with jumps}


\selectlanguage{english}
\author[authorlabel1]{Arturo Kohatsu-Higa},
\ead{khts00@fc.ritsumei.ac.jp}
\author[authorlabel2]{Eulalia Nualart},
\ead{eulalia@nualart.es}
\author[authorlabel3]{Ngoc Khue Tran}
\ead{trankhue@math.univ-paris13.fr}
\address[authorlabel1]{Department of Mathematical Sciences - Ritsumeikan University and Japan Science and Technology Agency, 1-1-1 Nojihigashi, Kusatsu, Shiga, 525-8577, Japan}
\address[authorlabel2]{Dept. Economics and Business, Universitat Pompeu Fabra and Barcelona Graduate School of Economics, Ram\'on Trias Fargas 25-27, 08005 Barcelona, Spain}
\address[authorlabel3]{Universit\'e Paris 13, Sorbonne Paris Cit\'e, LAGA, CNRS, (UMR 7539), F-93430 Villetaneuse, France}


\begin{abstract}
\selectlanguage{english}
In this paper, we consider a linear model with jumps driven by a Brownian motion and a compensated Poisson process, whose drift and diffusion coefficients as well as its intensity are unknown parameters. Supposing that the process is observed  discretely at high frequency we derive the local asymptotic normality (LAN) property. In order to obtain this result, Malliavin calculus and Girsanov's theorem are applied in order to write the log-likelihood ratio in terms of sums of conditional expectations, for which a central limit theorem for triangular arrays can be applied.

\vskip 0.5\baselineskip

\selectlanguage{francais}
\noindent{\bf R\'esum\'e} \vskip 0.5\baselineskip \noindent
{\bf La propri\'et\'e LAN pour un mod\`ele lin\'eaire avec sauts.} Dans cet article, nous consid\'erons un mod\`ele lin\'eaire avec sauts dirig\'e par un mouvement Brownien et un processus de Poisson compens\'e dont les coefficients et l'intensit\'e d\'ependent de param\`etres inconnus. Supposant que le processus est observ\'e \`a haute fr\'equence, nous obtenons la propri\'et\'e de normalit\'e asymptotique locale. Pour cela, le calcul de Malliavin et le th\'eor\`eme de Girsanov sont appliqu\'es afin d'\'ecrire le logarithme du rapport de vraisemblances comme une somme d'esp\'erances conditionnelles, pour laquelle un th\'eor\`eme centrale limite pour des suites triangulaires peut \^etre appliqu\'e.

\end{abstract}
\end{frontmatter}


\selectlanguage{english}
\section{Introduction and main result}
On a complete probability space $(\Omega, \mathcal{F}, \P)$, we consider the following stochastic process $X^{\theta,\sigma,\lambda}=(X_t^{\theta,\sigma,\lambda})_{t \geq 0}$ in $\R$ defined by
\begin{equation}\label{eq1}
X_t^{\theta,\sigma,\lambda}=x+\theta t +\sigma B_t + N_t -\lambda t,
\end{equation}
where $B=(B_t)_{t \geq 0}$ is a standard Brownian motion, $N=(N_t)_{t \geq 0}$ is a Poisson process with intensity $\lambda >0$ independent of $B$, and we denote by $(\widetilde{N}_{t}^{\lambda})_{t \geq 0}$ the compensated Poisson process  $\widetilde{N}_{t}^{\lambda}:=N_t-\lambda t$. The parameters $(\theta,\sigma,\lambda)\in \Theta\times \Sigma \times\Lambda$ are unknown and $\Theta, \Sigma$ and $\Lambda$ are closed intervals of $\R, \R^{\ast}_+$ and $\R^{\ast}_+$, where $\R^{\ast}_+=\R_+\setminus\{0\}$. Let $\{\mathcal{F}_t\}_{t\geq0}$ denote the natural filtration generated by $B$ and $N$. We denote by $\P_x^{\theta,\sigma,\lambda}$ the probability law induced by the $\mathcal{F}$-adapted c\`{a}dl\`{a}g stochastic process $X^{\theta,\sigma,\lambda}$ starting at $x$, and by $\E_x^{\theta,\sigma,\lambda}$ the expectation with respect to $\P_x^{\theta,\sigma,\lambda}$. Let $\overset{\P_x^{\theta,\sigma,\lambda}}{\longrightarrow}$ and $\overset{\mathcal{L}(\P_x^{\theta,\sigma,\lambda})}{\longrightarrow}$ denote the convergence in $\P_x^{\theta,\sigma,\lambda}$-probability and in $\P_x^{\theta,\sigma,\lambda}$-law, respectively.

For $(\theta,\sigma,\lambda)\in \Theta\times \Sigma \times\Lambda$, we consider an equidistant discrete observation of the process $X^{\theta,\sigma,\lambda}$ which is denoted by  $X^{n}=(X_{t_0}, X_{t_1},...,X_{t_n})$, where $t_k=k \Delta_n$ for $k \in \{0,...,n\}$, and $\Delta_n\leq 1$. We assume that the high-frequency observation condition holds. That is,
\begin{equation} \label{eq4}
n \Delta_n \rightarrow \infty, \quad \text{ and } \quad \Delta_n \rightarrow 0, \quad \text{as } n \rightarrow \infty.
\end{equation}
Let $p(\cdot;(\theta,\sigma,\lambda))$ denote the density of the random vector $X^{n}$ under the parameter $(\theta,\sigma,\lambda)$. For $(u,v,w)\in\R^3$, set $\theta_n:=\theta+\dfrac{u}{\sqrt{n\Delta_n}}, \sigma_n:=\sigma+\dfrac{v}{\sqrt{n}}, \lambda_n:=\lambda+\dfrac{w}{\sqrt{n\Delta_n}}$.

The aim of this paper is to prove the following LAN property.
\begin{theorem}\label{theorem}
Assume condition \eqref{eq4}. Then, the LAN property holds for all $(\theta,\sigma,\lambda)\in \Theta\times \Sigma \times\Lambda$ with rate of convergence $(\sqrt{n\Delta_n},\sqrt{n},\sqrt{n\Delta_n})$ and asymptotic Fisher information matrix $\Gamma(\theta,\sigma,\lambda)$. That is, for all $z=(u,v,w)\in\R^3$, as $n\to\infty$,
\begin{equation*} 
\log\frac{p\left(X^{n};\left(\theta_n,\sigma_n,\lambda_n\right)\right)}{p\left(X^{n};\left(\theta,\sigma,\lambda\right)\right)}\\
\overset{\mathcal{L}(\P_{x}^{\theta,\sigma,\lambda})}{\longrightarrow} z^\mathsf{T}\mathcal{N}\left(0,\Gamma(\theta,\sigma,\lambda)\right)-\dfrac{1}{2}z^\mathsf{T} \Gamma(\theta,\sigma,\lambda)z,
\end{equation*}
where $\mathcal{N}(0,\Gamma(\theta,\sigma,\lambda))$ is a centered $\R^3$-valued Gaussian vector with covariance matrix 
$$\Gamma(\theta,\sigma,\lambda)=\dfrac{1}{\sigma^2}\begin{pmatrix}1&0&-1\\0&2&0\\
-1&0&1+\frac{\sigma^2}{\lambda}
\end{pmatrix}.$$
\end{theorem}

Theorem \ref{theorem} extends in the linear case and in the presence of jumps the results of Gobet in \cite{G01} and \cite{G02} for multidimensional continuous elliptic diffusions. The main idea of these papers is to use the Malliavin calculus in order to obtain an expression for the derivative of the log-likelihood function in terms of a conditional expectation. Some extensions of Gobet's work with the presence of jumps are given for e.g. in \cite{CDG}, \cite{G3}, and \cite{K13}. However, in the present note, we estimate the coefficients and jump intensity parameters  at the same time. 
The main motivation for this paper is to show some of the important properties and arguments in order to prove the LAMN property in the non-linear case whose proof is non-trivial.  In particular, we present four important Lemmas of independent interest which will be key elements in dealing with the non-linear case. The key argument consists in conditioning on the number of jumps within the conditional expectation which expresses the transition density and outside it. When these two conditionings relate to different jumps one may use a large deviation principle in the estimate. When they are equal one uses the complementary set. Within all these arguments the Gaussian type upper and lower bounds of the density conditioned on the jumps is again strongly used. This idea seems to have many other uses in the set-up of stochastic differential equations driven by a Brownian motion and a jump process. We remark here that a plain It\^o-Taylor expansion would not solve the problem as higher moments of the Poisson process do not become smaller as the expansion order increases.

\section{Preliminaries}

In this Section we introduce the preliminary results needed for the proof of Theorem \ref{theorem}.
In order to deal with the likelihood ratio in Theorem \ref{theorem}, we will use the following decomposition 
\begin{equation}
\label{eq:dec} 
\log\frac{p\left(X^{n};\left(\theta_n,\sigma_n,\lambda_n\right)\right)}{p\left(X^{n};\left(\theta,\sigma,\lambda\right)\right)}=\log\dfrac{p\left(X^{n};\left(\theta_n,\sigma_n,\lambda_n\right)\right)}{p\left(X^{n};\left(\theta_n,\sigma,\lambda_n\right)\right)}+\log\dfrac{p\left(X^{n};\left(\theta_n,\sigma,\lambda_n\right)\right)}{p\left(X^{n};\left(\theta_n,\sigma,\lambda\right)\right)}+\log\dfrac{p\left(X^{n};\left(\theta_n,\sigma,\lambda\right)\right)}{p\left(X^{n};\left(\theta,\sigma,\lambda\right)\right)}.
\end{equation}

For each of the above terms we will use a mean value theorem and then analyze each term.
We start as in Gobet \cite{G01} applying the integration by parts formula of the Malliavin calculus on each interval $[t_k, t_{k+1}]$ to obtain the following expressions for the derivatives of the log-likelihood function w.r.t. $\theta$ and $\sigma$. Moreover, using Girsanov's theorem, we obtain the following expression for the  log-likelihood function w.r.t. $\lambda$. For any $t>s$, we denote by $p^{\theta,\sigma,\lambda}(t-s,x,y)$ the transition density of $X_t^{\theta,\sigma,\lambda}$ conditioned on $X_s^{\theta,\sigma,\lambda}=x$.

\begin{proposition} \label{prop1}
For all $\theta\in\R, \sigma, \lambda\in\R^{\ast}_+$, and $k \in \{0,...,n-1\}$,
\begin{equation*} \begin{split}
&\dfrac{\partial_{\theta}p^{\theta,\sigma,\lambda}}{p^{\theta,\sigma,\lambda}}\left(\Delta_n,X_{t_k},X_{t_{k+1}}\right)=\dfrac{1}{\sigma}\E_{X_{t_k}}^{\theta,\sigma,\lambda}\left[B_{t_{k+1}}-B_{t_{k}}\Big\vert X_{t_{k+1}}^{\theta,\sigma,\lambda}=X_{t_{k+1}}\right],\\
&\dfrac{\partial_{\sigma}p^{\theta,\sigma,\lambda}}{p^{\theta,\sigma,\lambda}}\left(\Delta_n,X_{t_k},X_{t_{k+1}}\right)=\dfrac{1}{\Delta_n}\E_{X_{t_k}}^{\theta,\sigma,\lambda}\left[\left(B_{t_{k+1}}-B_{t_{k}}\right)^2\Big\vert X_{t_{k+1}}^{\theta,\sigma,\lambda}=X_{t_{k+1}}\right] -\frac{1}{\sigma}, \\
&\dfrac{\partial_{\lambda}p^{\theta,\sigma,\lambda}}{p^{\theta,\sigma,\lambda}}(\Delta_n,X_{t_k},X_{t_{k+1}})=\E_{X_{t_k}}^{\theta,\sigma,\lambda} \left[-\dfrac{B_{t_{k+1}}-B_{t_{k}}}{\sigma}+\dfrac{\widetilde{N}_{t_{k+1}}^{\lambda}-\widetilde{N}_{t_k}^{\lambda}}{\lambda}\bigg\vert X_{t_{k+1}}^{\theta,\sigma,\lambda}=X_{t_{k+1}}\right].
\end{split}
\end{equation*}
\end{proposition}

We next present the four Lemmas mentionned in the Introduction.
Consider the events $J_m:=\{N_{t_{k+1}}-N_{t_{k}}=m\}$, for all $m\geq 0$ and $k \in \{0,...,n-1\}$. 
\begin{lemma} \label{lemma4}
For all $\theta\in\R, \sigma, \lambda\in\R^{\ast}_+,\ k \in \{0,...,n-1\}$, and $m\geq 0$,
\begin{equation*}\begin{split}
\P_{X_{t_k}}^{\theta,\sigma,\lambda}\left(J_m\Big\vert X_{t_{k+1}}^{\theta,\sigma,\lambda}=X_{t_{k+1}}\right)=\dfrac{ e^{-\left(X_{t_{k+1}}-X_{t_k}-m-\left(\theta-\lambda\right)\Delta_n\right)^2/(2\sigma^2\Delta_n)}\frac{(\lambda\Delta_n)^m}{m!}}{\sum_{i=0}^{\infty} e^{-\left(X_{t_{k+1}}-X_{t_k}-i-\left(\theta-\lambda\right)\Delta_n\right)^2/(2\sigma^2\Delta_n)}\frac{(\lambda\Delta_n)^i}{i!}}.
\end{split}
\end{equation*}
\end{lemma}

For all $j, p\geq 0$ and $k \in \{0,...,n-1\}$, we introduce the random variable
\begin{equation*} \begin{split}
S_{j}^{p}:={\bf 1}_{J_j}\E_{X_{t_k}}^{\bar{\theta},\bar{\sigma},\bar{\lambda}}\left[{\bf 1}_{J_j^c}\left(N_{t_{k+1}}-N_{t_{k}}\right)^p\Big\vert X_{t_{k+1}}^{\bar{\theta},\bar{\sigma},\bar{\lambda}}=X_{t_{k+1}}\right].
\end{split}
\end{equation*}
We remark that heuristically the indicator functions ${\bf 1}_{J_j}$ and ${\bf 1}_{J_j^c}$ outside and inside the conditional expectation correspond to restrictions on $X^{(\theta,\sigma,\lambda)}$ and $X^{(\bar{\theta},\bar{\sigma},\bar{\lambda})}$, respectively. 

\begin{lemma}\label{lemma5} For all $\theta, \bar{\theta}\in\R, \sigma, \bar{\sigma}, \lambda, \bar{\lambda}\in \R^{\ast}_+, j, p\geq 0$ and $k \in \{0,...,n-1\}$,
\begin{equation}\label{eq44} \begin{split}
S_{j}^{p}={\bf 1}_{J_j}\dfrac{\sum_{m=0: m\neq j}^{\infty}m^pe^{-\left(\sigma\left(B_{t_{k+1}}-B_{t_k}\right)+j-m+(\theta-\bar{\theta}-\lambda+\bar{\lambda})\Delta_n\right)^2/(2\bar{\sigma}^2\Delta_n)}\frac{(\bar{\lambda}\Delta_n)^m}{m!}}{\sum_{i=0}^{\infty}e^{-\left(\sigma\left(B_{t_{k+1}}-B_{t_k}\right)+j-i+(\theta-\bar{\theta}-\lambda+\bar{\lambda})\Delta_n\right)^2/(2\bar{\sigma}^2\Delta_n)}\frac{(\bar{\lambda}\Delta_n)^i}{i!}}.
\end{split}
\end{equation}
\end{lemma}

We next fix $\alpha\in(0,\frac{1}{2})$, and analyze $S^p_j$ in two separate cases as follows
\begin{equation*} 
S_j^p=S_j^p {\bf 1}_{\{\vert B_{t_{k+1}}-B_{t_k}\vert\leq \Delta_n^{\alpha}\}}+ S_j^p{\bf 1}_{\{\vert B_{t_{k+1}}-B_{t_k}\vert > \Delta_n^{\alpha}\}}=:S_{1,j}^{p}+S_{2,j}^{p}.
\end{equation*}
Furthermore, we write
$
S_{1,j}^{p}=S_{1,1,j}^{p}+S_{1,2,j}^{p},
$
and $
S_{2,j}^{p}=S_{2,1,j}^{p}+S_{2,2,j}^{p}$, 
where  $S_{1,1,j}^{p}$ and $S_{2,1,j}^{p}$ contain the terms $\sum_{m<j}$,
and $S_{1,2,j}^{p}$ and  $S_{2,2,j}^{p}$ contain the terms $\sum_{m>j}$ in (\ref{eq44}).
\begin{lemma}\label{lemma6} Assume that $\vert\theta-\bar{\theta}\vert\leq\frac{C}{\sqrt{n\Delta_n}}$ and $\vert\lambda-\bar{\lambda}\vert\leq\frac{C}{\sqrt{n\Delta_n}}$, for some constant $C>0$. Then for all $\sigma, \bar{\sigma}\in \R^{\ast}_+, j, p\geq 0,\ k \in \{0,...,n-1\}$, and for $n$ large enough,
\begin{align*} 
S_{1,1,j}^{p}&\leq {\bf 1}_{J_j}\dfrac{j!}{(\bar{\lambda}\Delta_n)^j}\sum_{m<j}m^pe^{-\frac{(j-m)^2}{4\bar{\sigma}^2\Delta_n}}\frac{(\bar{\lambda}\Delta_n)^m}{m!},\qquad 
S_{1,2,j}^{p}\leq {\bf 1}_{J_j}e^{-\frac{1}{4\bar{\sigma}^2\Delta_n}}\sum_{\ell>0}(\ell+j)^p\frac{(\bar{\lambda}\Delta_n)^{\ell}}{\ell!}, \\
S_{2,1,j}^{p}&\leq j^p{\bf 1}_{J_j}{\bf 1}_{\{\vert B_{t_{k+1}}-B_{t_k}\vert>\Delta_n^{\alpha}\}},\qquad
S_{2,2,j}^{p}\leq {\bf 1}_{J_j}{\bf 1}_{\{\vert B_{t_{k+1}}-B_{t_k}\vert>\Delta_n^{\alpha}\}}\sum_{\ell=0}^{\infty}(\ell+j+1)^p\frac{(\bar{\lambda}\Delta_n)^{\ell}}{\ell!}.
\end{align*}
\end{lemma}

For all $p\geq 0$ and $k \in \{0,...,n-1\}$, set
\begin{equation*} \begin{split}
M_{1,p}^{\bar{\theta},\bar{\sigma},\bar{\lambda}}:&=\sum_{j=0}^{\infty}j^p\E_{X_{t_k}}^{\theta,\sigma,\lambda}\left[{\bf 1}_{J_j}\E_{X_{t_k}}^{\bar{\theta},\bar{\sigma},\bar{\lambda}}\left[{\bf 1}_{J_j^c}\Big\vert X_{t_{k+1}}^{\bar{\theta},\bar{\sigma},\bar{\lambda}}=X_{t_{k+1}}\right]\right],\\
M_{2,p}^{\bar{\theta},\bar{\sigma},\bar{\lambda}}:&=\sum_{j=0}^{\infty}\E_{X_{t_k}}^{\theta,\sigma,\lambda}\left[{\bf 1}_{J_j}\E_{X_{t_k}}^{\bar{\theta},\bar{\sigma},\bar{\lambda}}\left[{\bf 1}_{J_j^c}\left(N_{t_{k+1}}-N_{t_{k}}\right)^p\Big\vert X_{t_{k+1}}^{\bar{\theta},\bar{\sigma},\bar{\lambda}}=X_{t_{k+1}}\right]\right].
\end{split}
\end{equation*}

\begin{lemma}\label{lemma7} Assume that $\vert\theta-\bar{\theta}\vert\leq\frac{C}{\sqrt{n\Delta_n}}$ and $\vert\lambda-\bar{\lambda}\vert\leq\frac{C}{\sqrt{n\Delta_n}}$, for some constant $C>0$. Then, for any $\sigma, \bar{\sigma}\in \Sigma$, $p\geq 0$,  and for $n$ large enough, there exist constants $C_1, C_2>0$ such that for all $\alpha\in(0,\frac{1}{2})$, and  $\ k \in \{0,...,n-1\}$, 
\begin{equation*} \begin{split}
M_{1,p}^{\bar{\theta},\bar{\sigma},\bar{\lambda}}+ M_{2,p}^{\bar{\theta},\bar{\sigma},\bar{\lambda}}\leq C_1e^{-\frac{1}{C_2\Delta_n^{1-2\alpha}}}.
\end{split}
\end{equation*}
\end{lemma}

We next recall a convergence in probability result, and a central limit theorem for triangular arrays of random variables. For each $n\in\mathbb{N}$, let $(Z_{k,n})_{k\geq 1}$ and $(\zeta_{k,n})_{k\geq 1}$ be two sequences of random variables defined on the filtered probability space $(\Omega, \mathcal{F}, (\mathcal{F}_t)_{t\geq 0}, \P)$, and assume that they are $\mathcal{F}_{t_{k+1}}$-measurable.
\begin{lemma}\label{zero} \textnormal{\cite[Lemma 9]{GJ93}} Assume that 
$\sum_{k=0}^{n-1}\E\left[Z_{k,n}\vert \mathcal{F}_{t_k}\right] \overset{\P}{\longrightarrow} 0$, and $\sum_{k=0}^{n-1}\E\left[Z_{k,n}^2\vert \mathcal{F}_{t_k} \right]\overset{\P}{\longrightarrow} 0$, as $n  \rightarrow \infty$.
Then 
$
\sum_{k=0}^{n-1}Z_{k,n}\overset{\P}{\longrightarrow} 0,
$ as $n  \rightarrow \infty$.
\end{lemma}

\begin{lemma}\label{clt} \textnormal{\cite[Lemma 4.3]{J11}}  
Assume that there exist real numbers $M$ and $V>0$ such that as $n  \rightarrow \infty$,  
\begin{equation*} \begin{split}
\sum_{k=0}^{n-1}\E\left[\zeta_{k,n}\vert \mathcal{F}_{t_k}\right] \overset{\P}{\longrightarrow} M, \quad \sum_{k=0}^{n-1}\left(\E\left[\zeta_{k,n}^2\vert \mathcal{F}_{t_k} \right]-\left(\E\left[\zeta_{k,n}\vert \mathcal{F}_{t_k}\right]\right)^2\right)\overset{\P}{\longrightarrow} V, \text{ and } 
\sum_{k=0}^{n-1}\E\left[\zeta_{k,n}^4\vert \mathcal{F}_{t_k}\right] \overset{\P}{\longrightarrow} 0.
\end{split}
\end{equation*}
Then as $n  \rightarrow \infty$,  $\sum_{k=0}^{n-1}\zeta_{k,n}\overset{\mathcal{L}(\P)}{\longrightarrow} \mathcal{N}+M$, where $\mathcal{N}$ is a centered  Gaussian variable with variance $V$.
\end{lemma}

\section{Proof of Theorem \ref{theorem}}

For $\ell\in[0,1]$, set $\theta(\ell):=\theta_n(\ell,u):=\theta+\dfrac{\ell u}{\sqrt{n\Delta_n}}, \sigma(\ell):=\sigma_n(\ell,v):=\sigma+\dfrac{\ell v}{\sqrt{n}}, \lambda(\ell):=\lambda_n(\ell,w):=\lambda+\dfrac{\ell w}{\sqrt{n\Delta_n}}$.
Applying the Markov property and Proposition \ref{prop1} to each term in (\ref{eq:dec}), we obtain that 
\begin{equation*}  \begin{split}
\log\dfrac{p\left(X^{n};\left(\theta_n,\sigma,\lambda\right)\right)}{p\left(X^{n};\left(\theta,\sigma,\lambda\right)\right)}&=\sum_{k=0}^{n-1}\log\dfrac{p^{\theta_n,\sigma,\lambda}}{p^{\theta,\sigma,\lambda}}\left(\Delta_n,X_{t_k},X_{t_{k+1}}\right)\\
&=\sum_{k=0}^{n-1}\dfrac{u}{\sqrt{n\Delta_n}}\int_0^1\dfrac{\partial_{\theta}p^{\theta(\ell),\sigma,\lambda}}{p^{\theta(\ell),\sigma,\lambda}}\left(\Delta_n,X_{t_k},X_{t_{k+1}}\right)d\ell \\
&=\sum_{k=0}^{n-1}\dfrac{u}{\sqrt{n\Delta_n}}\dfrac{1}{\sigma}\int_0^1\E_{X_{t_k}}^{\theta(\ell),\sigma,\lambda}\left[B_{t_{k+1}}-B_{t_{k}}\Big\vert X_{t_{k+1}}^{\theta(\ell),\sigma,\lambda}=X_{t_{k+1}}\right]d\ell,
\end{split}
\end{equation*}
\begin{equation*}  \begin{split}
&\log\dfrac{p\left(X^{n};\left(\theta_n,\sigma_n,\lambda_n\right)\right)}{p\left(X^{n};\left(\theta_n,\sigma,\lambda_n\right)\right)}=\sum_{k=0}^{n-1}\dfrac{v}{\sqrt{n}}\int_0^1\dfrac{\partial_{\sigma}p^{\theta_n,\sigma(\ell),\lambda_n}}{p^{\theta_n,\sigma(\ell),\lambda_n}}\left(\Delta_n,X_{t_k},X_{t_{k+1}}\right)d\ell\\
&=\sum_{k=0}^{n-1}\dfrac{v}{\sqrt{n}}\int_0^1 \left(\dfrac{1}{\Delta_n}\E_{X_{t_k}}^{\theta_n,\sigma(\ell),\lambda_n}\left[\left(B_{t_{k+1}}-B_{t_{k}}\right)^2\Big\vert X_{t_{k+1}}^{\theta_n,\sigma(\ell),\lambda_n}=X_{t_{k+1}}\right] -\frac{1}{\sigma(\ell)} \right)d\ell,
\end{split}
\end{equation*}
and
\begin{equation*}  \begin{split}
&\log\dfrac{p\left(X^{n};\left(\theta_n,\sigma,\lambda_n\right)\right)}{p\left(X^{n};\left(\theta_n,\sigma,\lambda\right)\right)}=\sum_{k=0}^{n-1}\dfrac{w}{\sqrt{n\Delta_n}}\int_0^1\dfrac{\partial_{\lambda}p^{\theta_n,\sigma,\lambda(\ell)}}{p^{\theta_n,\sigma,\lambda(\ell)}}\left(\Delta_n,X_{t_k},X_{t_{k+1}}\right)d\ell\\
&=\sum_{k=0}^{n-1}\dfrac{w}{\sqrt{n\Delta_n}}\int_0^1\E_{X_{t_k}}^{\theta_n,\sigma,\lambda(\ell)} \left[-\dfrac{B_{t_{k+1}}-B_{t_{k}}}{\sigma}+\dfrac{\widetilde{N}_{t_{k+1}}^{\lambda(\ell)}-\widetilde{N}_{t_k}^{\lambda(\ell)}}{\lambda(\ell)}\bigg\vert X_{t_{k+1}}^{\theta_n,\sigma,\lambda(\ell)}=X_{t_{k+1}}\right]d\ell.
\end{split}
\end{equation*}
Now using equation \eqref{eq1}, we obtain the following expansion of the log-likelihood ratio
\begin{equation*} 
\log\frac{p\left(X^{n};\left(\theta_n,\sigma_n,\lambda_n\right)\right)}{p\left(X^{n};\left(\theta,\sigma,\lambda\right)\right)}=\sum_{k=0}^{n-1}\left(\xi_{k,n}+H_{k,n}+\eta_{k,n}+M_{k,n}+\beta_{k,n}-R_{k,n} \right),
\end{equation*}
where
\begin{equation*}\begin{split}
\xi_{k,n}&:=\dfrac{u}{\sqrt{n\Delta_n}}\dfrac{1}{\sigma^2}\left(\sigma\left(B_{t_{k+1}}-B_{t_{k}}\right)-\dfrac{u\Delta_n}{2\sqrt{n\Delta_n}}\right),\\
H_{k,n}&:=\dfrac{u}{\sqrt{n\Delta_n}}\dfrac{1}{\sigma^2}\left(\widetilde{N}_{t_{k+1}}^{\lambda}-\widetilde{N}_{t_{k}}^{\lambda}-\int_0^1\E_{X_{t_k}}^{\theta(\ell),\sigma,\lambda}\left[\widetilde{N}_{t_{k+1}}^{\lambda}-\widetilde{N}_{t_{k}}^{\lambda}\Big\vert X_{t_{k+1}}^{\theta(\ell),\sigma,\lambda}=X_{t_{k+1}}\right]d\ell\right), \\
\eta_{k,n}&:=\dfrac{v}{\sqrt{n}}\int_0^1\dfrac{1}{\Delta_n}\left(\dfrac{\sigma^2}{\sigma(\ell)^3}\left(B_{t_{k+1}}-B_{t_{k}}\right)^2-\dfrac{\Delta_n}{\sigma(\ell)}\right)d\ell,\\
M_{k,n}&:=\dfrac{v}{\sqrt{n}}\int_0^1\dfrac{1}{\Delta_n}\dfrac{1}{\sigma(\ell)^3}\bigg\{\left(\theta\Delta_n+\widetilde{N}_{t_{k+1}}^{\lambda}-\widetilde{N}_{t_{k}}^{\lambda}\right)^2+2\sigma\left(B_{t_{k+1}}-B_{t_{k}}\right)\left(\theta\Delta_n+\widetilde{N}_{t_{k+1}}^{\lambda}-\widetilde{N}_{t_{k}}^{\lambda}\right)\\
&\qquad\qquad-\E_{X_{t_k}}^{\theta_n,\sigma(\ell),\lambda_n} \bigg[\left(\theta_n\Delta_n+\widetilde{N}_{t_{k+1}}^{\lambda_n}-\widetilde{N}_{t_{k}}^{\lambda_n}\right)^2\\
&\qquad\qquad+2\sigma(\ell)\left(B_{t_{k+1}}-B_{t_{k}}\right)\left(\theta_n\Delta_n+\widetilde{N}_{t_{k+1}}^{\lambda_n}-\widetilde{N}_{t_{k}}^{\lambda_n}\right) \bigg\vert X_{t_{k+1}}^{\theta_n,\sigma(\ell),\lambda_n}=X_{t_{k+1}}\bigg]\bigg\}d\ell, \\
\beta_{k,n}&:=-\dfrac{w}{\sqrt{n\Delta_n}}\dfrac{1}{\sigma^2}\left(\sigma\left(B_{t_{k+1}}-B_{t_{k}}\right)+\dfrac{w\Delta_n}{2\sqrt{n\Delta_n}}-\dfrac{u\Delta_n}{\sqrt{n\Delta_n}}\right)\\
&\qquad+\dfrac{w}{\sqrt{n\Delta_n}}\int_0^1\E_{X_{t_k}}^{\theta_n,\sigma,\lambda(\ell)}\left[\dfrac{\widetilde{N}_{t_{k+1}}^{\lambda(\ell)}-\widetilde{N}_{t_k}^{\lambda(\ell)}}{\lambda(\ell)}\bigg\vert X_{t_{k+1}}^{\theta_n,\sigma,\lambda(\ell)}=X_{t_{k+1}}\right]d\ell,\\
R_{k,n}&:=\dfrac{w}{\sqrt{n\Delta_n}}\dfrac{1}{\sigma^2}\int_0^1\left(\widetilde{N}_{t_{k+1}}^{\lambda(\ell)}-\widetilde{N}_{t_k}^{\lambda(\ell)} -\E_{X_{t_k}}^{\theta_n,\sigma,\lambda(\ell)}\left[  \widetilde{N}_{t_{k+1}}^{\lambda(\ell)}-\widetilde{N}_{t_k}^{\lambda(\ell)} \bigg\vert X^{\theta_n,\sigma,\lambda(\ell)}_{t_{k+1}}=X_{t_{k+1}}\right]\right)d\ell.
\end{split}
\end{equation*}

We next show that the random variables $\xi_{k,n}, \eta_{k,n}, \beta_{k,n}$ are the terms that contribute to the limit in Theorem \ref{theorem}, and  $H_{k,n}, M_{k,n}$ and $R_{k,n}$ are the negligible contributions.
Indeed, using Girsanov's theorem and Lemma \ref{lemma7}, we can show that the conditions of Lemma \ref{zero} under $\P_x^{\theta,\sigma,\lambda}$ hold for each term $H_{k,n}, M_{k,n}$ and $R_{k,n}$. That is,
\begin{lemma}\label{lemma1}
Assume condition \eqref{eq4}. Then, as $n\to\infty$, 
\begin{align*}
 &\sum_{k=0}^{n-1} \left( H_{k,n} +M_{k,n} -R_{k,n}\right) \overset{\P_x^{\theta,\sigma,\lambda}}{\longrightarrow} 0 \\
&\sum_{k=0}^{n-1} \E^{\theta,\sigma,\lambda}\left[\xi_{k,n}+ \eta_{k,n}+\beta_{k,n}\vert \mathcal{F}_{t_k}\right]
\overset{\P_x^{\theta,\sigma,\lambda}}{\longrightarrow}-\dfrac{u^2}{2\sigma^2}-\dfrac{v^2}{2}\dfrac{2}{\sigma^2}-\dfrac{w^2}{2\sigma^2}\left(1+\dfrac{\sigma^2}{\lambda}\right)+\dfrac{uw}{\sigma^2}\\
&\sum_{k=0}^{n-1}\left(\E^{\theta,\sigma,\lambda}\left[\xi_{k,n}^2+\eta_{k,n}^2+\beta_{k,n}^2\vert \mathcal{F}_{t_k}\right]-\E^{\theta,\sigma,\lambda}\left[\xi_{k,n}\vert\mathcal{F}_{t_k}\right]^2-\E^{\theta,\sigma,\lambda}\left[\eta_{k,n}\vert \mathcal{F}_{t_k}\right]^2-\E^{\theta,\sigma,\lambda}\left[\beta_{k,n}\vert \mathcal{F}_{t_k}\right]^2\right) \\
&\qquad \qquad\qquad \qquad\overset{\P_x^{\theta,\sigma,\lambda}}{\longrightarrow}\dfrac{u^2}{\sigma^2}+2\dfrac{v^2}{\sigma^2}+\dfrac{w^2}{\sigma^2}\left(1+\dfrac{\sigma^2}{\lambda}\right)\\
&\sum_{k=0}^{n-1}
\E^{\theta,\sigma,\lambda}\left[\xi_{k,n}^4+\eta_{k,n}^4+\beta_{k,n}^4\vert \mathcal{F}_{t_k}\right]\overset{\P_x^{\theta,\sigma,\lambda}}{\longrightarrow}0 \\
&\sum_{k=0}^{n-1}\left(\E^{\theta,\sigma,\lambda}\left[\xi_{k,n}\eta_{k,n}\vert \mathcal{F}_{t_k}\right]-\E^{\theta,\sigma,\lambda}\left[\xi_{k,n}\vert \mathcal{F}_{t_k}\right]\E^{\theta,\sigma,\lambda}\left[\eta_{k,n}\vert \mathcal{F}_{t_k}\right]\right)\overset{\P_x^{\theta,\sigma,\lambda}}{\longrightarrow}0\\
&\sum_{k=0}^{n-1}\left(\E^{\theta,\sigma,\lambda}\left[\xi_{k,n}\beta_{k,n}\vert \mathcal{F}_{t_k}\right]-\E^{\theta,\sigma,\lambda}\left[\xi_{k,n}\vert \mathcal{F}_{t_k}\right]\E^{\theta,\sigma,\lambda}\left[\beta_{k,n}\vert \mathcal{F}_{t_k}\right]\right)\overset{\P_x^{\theta,\sigma,\lambda}}{\longrightarrow}-\dfrac{uw}{\sigma^2}\\
&\sum_{k=0}^{n-1}\left(\E^{\theta,\sigma,\lambda}\left[\eta_{k,n}\beta_{k,n}\vert \mathcal{F}_{t_k}\right]-\E^{\theta,\sigma,\lambda}\left[\eta_{k,n}\vert \mathcal{F}_{t_k}\right]\E^{\theta,\sigma,\lambda}\left[\beta_{k,n}\vert \mathcal{F}_{t_k}\right]\right)\overset{\P_x^{\theta,\sigma,\lambda}}{\longrightarrow}0. 
\end{align*}
\end{lemma}

Finally, Lemma \ref{clt} applied to  $\zeta_{k,n}=\xi_{k,n}+\eta_{k,n}+\beta_{k,n}$ concludes the proof of Theorem \ref{theorem}.

\section*{Acknowledgements}
Second  author acknowledges support from the European Union programme FP7-PEOPLE-2012-CIG under grant agreement 333938. Third author acknowledges support from the LIA CNRS Formath Vietnam and the program ARCUS MAE/IDF Vietnam.

\end{document}